\documentclass[reqno]{amsart}  
\theoremstyle{plain}

\newtheorem{theorem}{Theorem}[section]

\newtheorem{remark}[theorem]{Remark}
\newtheorem{proposition}[theorem]{Proposition}

\numberwithin{equation}{section}
\allowdisplaybreaks[1]

\theoremstyle{definition}

\theoremstyle{remark}

\newcommand{\cB}{{\mathcal B}}

\newcommand{\cD}{{\mathcal D}}
\newcommand{\cF}{{\mathcal F}}

\newcommand{\cL}{{\mathcal L}}

\newcommand{\cQ}{{\mathcal Q}}

\newcommand{\cU}{{\mathcal U}}
\newcommand{\cX}{{\mathcal X}}
\newcommand{\cY}{{\mathcal Y}}

\newcommand{\C}{{\mathbb C}}
\newcommand{\B}{{\mathbb B}}

\newcommand{\sbm}[1]{\left[\begin{smallmatrix} #1
         \end{smallmatrix}\right]}

\newcommand{\mat}[2]{\ensuremath{\left[\begin{array}{#1}
#2
\end{array} \right]}}

\begin{document}

\title[Equivalence of stability and performance]{
Equivalence of robust stabilization and robust performance via
feedback}
\author[J.A.~Ball]{Joseph A. Ball}
\address{Department of Mathematics,
Virginia Tech,
Blacksburg, VA 24061-0123, USA}
\email{ball@math.vt.edu}

\author[Q.~Fang]{Quanlei Fang}
\address{Department of Mathematics,
Virginia Tech,
Blacksburg, VA 24061-0123, USA}
\email{qlfang@math.vt.edu}

\author[G.~Groenewald]{Gilbert J.~Groenewald}
\address{Department of Mathematics, North West University,
Potchefstroom, 2520 South Africa}
\email{Gilbert.Groenewald@nwu.ac.za}

\author[S.~ter Horst]{Sanne Ter Horst}
\address{Department of Mathematics, Virginia Tech, Blacksburg, VA
24061-0123, USA}
\email{terhorst@math.vt.edu}



\begin{abstract}
One approach to robust control for linear plants with structured
uncertainty as well as for linear parameter-varying (LPV) plants
(where the controller has on-line access to the varying plant
parameters) is through linear-fractional-transformation (LFT) models.
Control issues to be addressed by controller design in this formalism
include robust stability and robust performance.  Here robust
performance is defined as the achievement of a uniform specified $L^{2}$-gain
tolerance for a disturbance-to-error map combined with robust
stability.  By setting the disturbance and error channels equal to
zero, it is clear that any criterion for robust performance also
produces a criterion for robust stability.  Counter-intuitively, as a
consequence of the so-called Main Loop Theorem, application of a
result on robust stability to a feedback configuration with an
artificial full-block uncertainty operator added in feedback
connection between the error and disturbance signals produces a
result on robust performance.  The main result here is that this
performance-to-stabilization
reduction principle must be handled with care for the case of dynamic feedback
compensation:  casual application of this principle leads to the
solution of a physically uninteresting problem, where the controller
is assumed to have access to the states in the artificially-added
feedback loop.  Application of the principle using a known more refined dynamic-control
robust stability criterion, where the user is allowed to specify controller
partial-state dimensions, leads to correct robust-performance results.
These latter results involve rank conditions in addition to Linear Matrix
Inequality (LMI) conditions.
\end{abstract}

\maketitle

\section{Introduction}
\label{sec:intro}

Linear-Fractional-Transformation (LFT) models have been used for the study of
stability issues
for systems with structured uncertainty
\cite{LZD, BeckDoyle},  of robust gain-scheduling for Linear
Parameter-Varying (LPV) systems \cite{Packard, LZD, AG, GA}, and of
model reduction for systems having
structured uncertainty \cite{BeckDoyleGlover, Beck, LiPaganini, BeckSCL06}.
It turns out that the LMI solution of the $H^{\infty}$-control problem
generalizes nicely to these more general structures; we refer the
reader to the books \cite{GN, DP} for nice expositions of these
and other related developments.

The results in the paper \cite{LZD} (see also \cite{LZDProc}) focus
on synthesis of controllers implementing a somewhat stronger notion
of stability known as $\cQ$-stability.  The notion of
$\cQ$-stability implies robust stability but the converse holds only
for special structures (see \cite{LZD}).  One such special structure
is the case where one allows the structured uncertainty to be
time-varying (and perhaps also causal and/or slowly time-varying in
a precise sense---see \cite{Paganini, BGM3}); then $\cQ$-stability
is equivalent to robust stability with respect to this enlarged
uncertainty structure. This observation gives perhaps the most
compelling system-theoretic interpretation of $Q$-stability.  Even
when one is not working with this enlarged uncertainty structure,
$\cQ$-stability is still attractive since it is sufficient for
robust stability and can be characterized in LMI form.

One result in \cite{LZD} is a characterization of the existence of a
static output feedback controller implementing $\cQ$-stability in
terms of the existence of positive-definite solutions $X,Y$ to a
pair of LMIs; the additional coupling condition $Y = X^{-1}$
destroys the convex character of the solution criterion and thereby
makes the solution criterion computationally unattractive.  A second
result provides an LMI characterization for the existence of a
dynamic (in the sense of multidimensional linear systems) controller
and provides a Youla parametrization for the set of all such
controllers.  The question of the existence of controllers for
LFT-model systems achieving $\cQ$-performance (a scaled version of
robust performance) is settled in \cite{Packard, GA, AG} (see the
book \cite{DP} for a nice overview); the existence of such
controllers is characterized in terms of the existence of structured
solutions $X,Y$ to a pair of LMIs subject to an additional coupling
condition
\begin{equation}  \label{coupling-Intro}
\begin{bmatrix} X & I \\ I & Y \end{bmatrix} \ge 0.
\end{equation}
Moreover, the rank of the various components of the controller
state-space can be prescribed by imposing additional rank
conditions on $\sbm{X & I \\ I & Y }$.

The purpose of this paper is to explain the precise logical
connections between results on robust stabilization versus results on
robust performance.  One direction is straightforward:  {\em any result on
robust performance gives rise to a result on robust stabilization by
specializing the robust performance result to the case where the
disturbance and error channels are trivial}. To recover the precise form of
the already existing results on robust stabilization however often
requires some additional algebraic manipulation.  The converse
direction is less obvious:  {\em any result on robust stabilization
implies a result on robust performance.} In its simplest form, as
pointed out in \cite{LZD}, this is
a consequence of the Main Loop Theorem for linear-fractional
maps (see \cite[Theorem 11.7, page 284]{ZDG}):

\medskip
\noindent
\textbf{Principle of Reduction of Robust Performance to Robust Stabilization:}
{\em robust performance can
be reduced to robust stability by adding a (fictitious) full-block uncertainty
feedback connection from the error channel to the disturbance channel}.

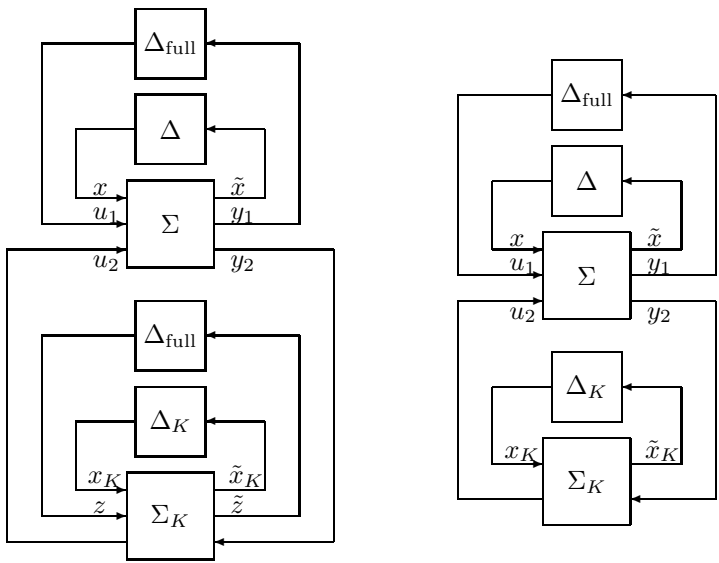
\begin{figure}[h]
\setlength{\unitlength}{0.09in}  
\hspace{2cm}
\begin{picture}(21,33)
\put(2,17){\framebox(5,5){$\Sigma$}}
\put(2.5,23){\framebox(4,4){$\Delta$}}
\put(2.5,28){\framebox(4,4){$\Delta_{\text{full}}$}}
\put(2,0){\framebox(5,5){$\Sigma_K$}}
\put(2.5,6){\framebox(4,4){$\Delta_K$}}
\put(2.5,11){\framebox(4,4){$\Delta_{\text{full}}$}}
\put(8,21.2){$\tilde x$}
\put(7,21){\line(1,0){3}}
\put(10,21){\line(0,1){4}}
\put(10,25){\vector(-1,0){3.5}}
\put(0,21.2){$x$}
\put(2.5,25){\line(-1,0){3.5}}
\put(-1,25){\line(0,-1){4}}
\put(-1,21){\vector(1,0){3}}
\put(8,19.8){$y_1$}
\put(7,19.5){\line(1,0){5}}
\put(12,19.5){\line(0,1){10.5}}
\put(12,30){\vector(-1,0){5.5}}
\put(0,19.8){$u_1$}
\put(2.5,30){\line(-1,0){5.5}}
\put(-3,30){\line(0,-1){10.5}}
\put(-3,19.5){\vector(1,0){5}}
\put(8,17){$y_2$}
\put(7,18){\line(1,0){7}}
\put(14,18){\line(0,-1){17}}
\put(14,1){\vector(-1,0){7}}
\put(0,17){$u_2$}
\put(2,1){\line(-1,0){7}}
\put(-5,1){\line(0,1){17}}
\put(-5,18){\vector(1,0){7}}
\put(7.9,4.4){$\tilde x_K$}
\put(7,4){\line(1,0){3}}
\put(10,4){\line(0,1){4}}
\put(10,8){\vector(-1,0){3.5}}
\put(-0.3,4.4){$x_K$}
\put(2.5,8){\line(-1,0){3.5}}
\put(-1,8){\line(0,-1){4}}
\put(-1,4){\vector(1,0){3}}
\put(8,2.7){$\tilde z$}
\put(7,2.5){\line(1,0){5}}
\put(12,2.5){\line(0,1){10.5}}
\put(12,13){\vector(-1,0){5.5}}
\put(0,2.7){$z$}
\put(2.5,13){\line(-1,0){5.5}}
\put(-3,13){\line(0,-1){10.5}}
\put(-3,2.5){\vector(1,0){5}}
\end{picture}
\hspace{.5cm}
\begin{picture}(19,33)
\put(2,14){\framebox(5,5){$\Sigma$}}
\put(2.5,20){\framebox(4,4){$\Delta$}}
\put(2.5,25){\framebox(4,4){$\Delta_{\text{full}}$}}
\put(2,2){\framebox(5,5){$\Sigma_K$}}
\put(2.5,8){\framebox(4,4){$\Delta_K$}}
\put(8,18.2){$\tilde x$}
\put(7,18){\line(1,0){3}}
\put(10,18){\line(0,1){4}}
\put(10,22){\vector(-1,0){3.5}}
\put(0,18.2){$x$}
\put(2.5,22){\line(-1,0){3.5}}
\put(-1,22){\line(0,-1){4}}
\put(-1,18){\vector(1,0){3}}
\put(8,16.8){$y_1$}
\put(7,16.5){\line(1,0){5}}
\put(12,16.5){\line(0,1){10.5}}
\put(12,27){\vector(-1,0){5.5}}
\put(0,16.8){$u_1$}
\put(2.5,27){\line(-1,0){5.5}}
\put(-3,27){\line(0,-1){10.5}}
\put(-3,16.5){\vector(1,0){5}}
\put(8,14){$y_2$}
\put(7,15){\line(1,0){5}}
\put(12,15){\line(0,-1){11.5}}
\put(12,3.5){\vector(-1,0){5}}
\put(0,14){$u_2$}
\put(2,3.5){\line(-1,0){5}}
\put(-3,3.5){\line(0,1){11.5}}
\put(-3,15){\vector(1,0){5}}
\put(7.9,5.9){$\tilde x_K$}
\put(7,5.5){\line(1,0){3}}
\put(10,5.5){\line(0,1){4.5}}
\put(10,10){\vector(-1,0){3.5}}
\put(-0.3,5.9){$x_K$}
\put(2.5,10){\line(-1,0){3.5}}
\put(-1,10){\line(0,-1){4.5}}
\put(-1,5.5){\vector(1,0){3}}
\end{picture}
\caption{Controller with and without access to the full block}
\label{fig:fullblock}    
\end{figure}

\medskip
\noindent The main point of the present paper is that this reduction
of robust performance to robust stability is not explained precisely
in the literature for the case of dynamic feedback for
multidimensional systems.  If one casually applies this principle to
the result from \cite{LZD} for the dynamic-feedback case, one
arrives at the results from \cite{Packard, AG} for robust
performance, but without the additional coupling constraint
\eqref{coupling-Intro}.  The explanation is that the condition with
the coupling constraint dropped does solve a robust-performance
problem, which, however, is a contrived problem of no physical
interest, namely the feedback configuration as on the left side of
Figure \ref{fig:fullblock}:~an LFT model $\Sigma$ for structured
uncertainty $\Delta$ with a controller $\Sigma_K$ that, besides the
controller-structured uncertainty $\Delta_K$, is granted access to
the artificial full-block uncertainty $\Delta_\textup{full}$ that
connects the error channel $y_1$ with the disturbance channel $u_1$.
Instead one must insist that the controller partial-state dimensions
for the states corresponding to the artificial full-block
uncertainty are zero so that the feedback configuration is as on the
right side of Figure \ref{fig:fullblock}. This additional constraint
on the dimension of the associated block of the controller state
space leads to the missing coupling condition. In this way the
reduction of robust performance to robust stabilization does hold,
but with proper attention paid to the controller information
structure.  Clarification of this point is the main contribution of
the present paper.

We mention that  another approach to robust control is to design a
controller to guarantee a uniform bound on the $L^{\infty}$-gain of
the disturbance-to-error map ($L^{1}$-control) rather than a uniform bound
on the $L^{2}$-gain of the disturbance-to-error map
($H^{\infty}$-control); a good overview for $L^{1}$-control is the
book \cite{D-DB}.  However, our focus here is on $H^{\infty}$-control
with the added feature that the disturbance/uncertainty is assumed to
have a structured form as given by an LFT model.

The paper is organized as follows. Section \ref{sec:intro} is the
present introduction. In Section \ref{sec:LFT} we review some known
results concerning LFT-model systems and  $\cQ$-stability and
$\cQ$-performance via output feedback. In particular, we recall a
theorem from \cite{DP} on $\cQ$-performance via output feedback with
controller partial-state dimension bounds. In the third section we
observe how this $\cQ$-performance result can be applied to obtain a
result on $\cQ$-stabilizability with dimension bounds on the controller;
in the two extreme cases where either (1) one demands that the
controller be static or (2) one imposes no restrictions on the size
of the partial states of the controller, auxiliary coupling conditions
in the criterion can be eliminated and we recover two
$\cQ$-stabilizability results from \cite{LZD} as corollaries. In Section
\ref{sec:Qperformance} we show that our theorem on $\cQ$-stabilizability
is actually equivalent to the $\cQ$-performance result in Section
\ref{sec:LFT}. We conclude this paper with a section on applications
for systems with structured uncertainty and for LPV systems.

\section{LFT-model systems in general}
\label{sec:LFT}

Let ${\mathbb F}$ be a field, taken to be either the complex numbers
${\mathbb C}$ or the real numbers ${\mathbb R}$. We define an {\em
LFT model for structured uncertainty} as follows. Assume that the
state space $\cX$, the input space $\cU$ and the output space $\cY$
are all finite-dimensional vector spaces over ${\mathbb F}$, say
 $$
 \cX = {\mathbb F}^{Z}, \quad \cU = {\mathbb F}^{M},\quad \cY={\mathbb F}^{N}.
 $$
 We then specify a direct-sum decomposition for $\cX$
   \begin{equation}  \label{cX}
    \cX = \oplus_{k=1}^{d} \cX_{k} \text{ with } \cX_{k} = {\mathbb
    F}^{n_{k} \cdot m_{k}},
    \, n_{1}\cdot m_{1} + \cdots + n_{d} \cdot m_{d} = Z.
   \end{equation}
   with associated {\em uncertainty structure} to be the collection
   $\boldsymbol \Delta$ of matrices of the block-diagonal form
   \begin{equation}  \label{scalar-blocks}
       \boldsymbol \Delta = \left\{
       \Delta = \begin{bmatrix} \Delta_{1} & & \\ & \ddots &
       \\ & & \Delta_{d}  \end{bmatrix} \colon
       \Delta_{k} = \begin{bmatrix} \delta_{k,11} I_{n_{k}} & \cdots &
       \delta_{k,1 m_{k}} I_{n_{k}} \\ \vdots & & \vdots \\
       \delta_{k,m_{k}1} I_{n_{k}} & \cdots & \delta_{k, m_{k} m_{k}}
       I_{n_{k}} \end{bmatrix} \right\},
       \end{equation}
       where $\delta_{k, ij}$ are arbitrary complex numbers for
       $k = 1, \dots, d$ and $i,j = 1, \dots,  m_{k}$.
       For short let us use the abbreviation
     $$
     \begin{bmatrix} \delta_{k,11} I_{n_{k}} & \cdots &
     \delta_{k,1 m_{k}} I_{n_{k}} \\ \vdots & & \vdots \\
     \delta_{k,m_{k}1} I_{n_{k}} & \cdots & \delta_{k, m_{k} m_{k}}
     I_{n_{k}} \end{bmatrix} =:
     \Delta^{0}_{k} \otimes  I_{n_{k}}
    $$
    where we have introduced the $m_{k} \times m_{k}$ matrix
    $\Delta^{0}_{k}$ with scalar entries given by
    \begin{equation}  \label{tensor}
    \Delta^{0}_{k} = \begin{bmatrix} \delta_{k,11} & \cdots &
    \delta_{k,1 m_{k}} \\ \vdots & & \vdots \\ \delta_{k,m_{k}1} &
    \cdots & \delta_{k,m_{k}m_{k}} \end{bmatrix}.
    \end{equation}
Then we define an {\em LFT model} (for structured uncertainty) to be any collection
of the form
$$\left\{ \begin{bmatrix} A & B \\ C & D \end{bmatrix} \colon
\begin{bmatrix} \cX \\ \cU \end{bmatrix} \to \begin{bmatrix} \cX \\
    \cY \end{bmatrix}, \quad {\boldsymbol \Delta}   \right\}.
$$
We shall also have use of the commutant of $\boldsymbol
\Delta$, denoted as $\cD_{\boldsymbol \Delta}$:
$$
   \cD_{\boldsymbol \Delta} =
\{ X \in \cL(\cX) \colon X \Delta = \Delta X \text{ for all }
   \Delta \in \boldsymbol \Delta \}.
$$
Explicitly, one can show that
the commutant $\cD_{\boldsymbol \Delta}$  consists of matrices $Q$ of the
$d \times d$-block diagonal form
  \begin{equation} \label{Dcl1}
    Q = \begin{bmatrix}
   Q_{1} & & \\ & \ddots & \\ & & Q_{d} \end{bmatrix}
   \end{equation}
   where, for $k = 1, \dots, d$,
   the $k$-th diagonal entry $Q_{k}$ in turn  has the  $m_{k} \times
   m_{k}$-block repeated diagonal form
   \begin{equation}  \label{Dcl2}
   Q_{k} = \begin{bmatrix} Q_{k,0} & & \\ & \ddots & \\ & &
   Q_{k,0} \end{bmatrix}
   \end{equation}
   where, finally,  the repeated block $Q_{k,0}$ is an arbitrary matrix of size  $n_{k} \times n_{k}$ with
   scalar entries.

 The associated {\em transfer function} in this context is the
   function of $\Delta \in \boldsymbol \Delta$ (defined at least for
   $\Delta$ having sufficiently small norm) given as the associated
   upper linear fractional transformation with symbol $\sbm{ A & B \\ C & D }$
   and load $\Delta$:
   \begin{equation}  \label{LFT-transfunc}
   \cF_{u}\left( \begin{bmatrix} A & B \\ C &  D \end{bmatrix}, \Delta \right)
   = D + C (I - \Delta A)^{-1} \Delta B \in \cL(\cU, \cY).
   \end{equation}
   Occasionally we shall also have use for the associated lower linear
fractional
   transformation with symbol $\sbm{ A & B \\ C & D }$ and load $\Delta'$:
   $$
   \cF_{\ell}\left( \begin{bmatrix} A & B \\ C &  D \end{bmatrix},
\Delta' \right)
   = A + B (I - \Delta' D)^{-1} \Delta' C \in \cL(\cX).
   $$

   This abstract notion of LFT model is used in \cite{LZD} (see the references
   there for more background) to model linear input/state/output
   systems having structured uncertainty.

   The classical case
   corresponds to the case ${\boldsymbol \Delta} = \{ \lambda
   I_{\cX}\}$ with $\lambda \in {\mathbb C}$ where $\lambda$ is the
   frequency variable; in this case, the upper linear fractional transformation
   is the transfer function of the
   discrete-time input/state/output linear system
   \begin{equation}  \label{sys}
   \left\{ \begin{array}{ccc}
   x(n+1) & = & A x(n) + B u(n) \\
     y(n) & = & C x(n) + D u(n)
     \end{array}  \right.\quad n=0,1,2,\ldots
   \end{equation}
   in the sense that $\cF_{u}\left( \sbm{A & B \\ C & D}, \lambda
   I_{\cX} \right)$ is the $Z$-transform
   $$
   \cF_{u}\left( \sbm{A & B \\ C & D}, \lambda
      I_{\cX}\right) = D + \sum_{n=1}^{\infty} CA^{n-1}B \lambda^{n}
    $$
   of the impulse response $\{ D, CB, CAB, \dots, CA^{n-1}B, \dots \}$
   of the system \eqref{sys},
   i.e., the output generated from zero initial condition and input
   signal corresponding to the unit impulse at time $0$ ($u(0) =
   I_{\cU}, u(n) = 0 $ for $n>0$).
   In case $\boldsymbol \Delta$ has the form
   \eqref{scalar-blocks} with ${\mathbb F}={\mathbb C}$ and $n_{k} = 1$ for all $k$, then
   $\cF_{u}\left( \sbm{A & B \\ C & D},  \Delta \right)$  can
   similarly be interpreted as the transfer function of a
   multidimensional linear system of Givone-Roesser type evolving on
   the integer lattice: here the frequency variable $\delta =
   (\delta_{1}, \dots, \delta_{d}) \in
   {\mathbb C}^{d}$ is $d$-dimensional (see \cite{GR}).  As an
   alternative interpretation, one can consider $\lambda:=
   \delta_{1}$ as the frequency variable for a 1-D system and the
   remaining parameters $\delta_{2}, \dots, \delta_{d}$ as values of
   parameters specifying a particular choice of disturbance within an
   admissible set of uncertainties.  Then $\cF_{u}
   \left( \sbm{A & B \\ C & D},  \Delta \right) $, considered as a
   function of $\lambda = \delta_{1}$ with the other $\delta$-values
   $\delta_{2}, \dots, \delta_{d}$ held fixed, specifies the classical
   transfer function for the system if one assumes the particular
   choice of uncertainty associated with the given fixed parameter
   values $\delta_{2}, \dots, \delta_{d}$ (see \cite{Packard,AG,DP}).
   One can even let $\delta_{1} \dots, \delta_{d}$ be formal noncommuting
   indeterminates and make
   sense of $\cF_{u}\left( \sbm{A & B \\ C & D}, \Delta \right)$ as a
   formal power series with coefficients equal to operators from $\cU$
   to $\cY$;
   then $\cF_{u}\left( \sbm{A & B \\ C & D}, \Delta \right)$ can be
   viewed as the transfer function of an input/state/output linear
   system having evolution along a free semigroup \cite{BGM1}.
   Alternatively, one can view the diagonal entries $\delta_{1},
   \dots, \delta_{d}$ as operators on $\ell^{2}$ with $\delta_{1}$
   equal to the shift operator, interpret $\delta_{2}, \dots,
   \delta_{d}$ as parameters associated with a particular choice of
   admissible structured time-varying disturbance in the system, and
   view the value of $\cF_{u}(\sbm{A & B \\ C & D}, \Delta)$ as
   the input-output map from $\cU \otimes \ell^{2}$ to $\cY \otimes \ell^{2}$
   for the system with
   particular choice of disturbance specified by  the choice of
   $\delta_{2}, \dots, \delta_{d}$.
   We discuss some of these various interpretations and their applications in more
   detail in our final Section \ref{sec:applications}.

   We now recall from \cite{LZD} how to formulate  robust
   stability and robust performance, along with the related notions
   of $\cQ$-stability and $\cQ$-performance, in the general context
   of an LFT model. Given an LFT model $\Sigma = ( \sbm{A & B \\ C & D},
\boldsymbol \Delta)$ we make the following definitions:
\begin{enumerate}
       \item The LFT model $\Sigma$  is {\em robustly stable} if
       $(I - \Delta A)$ is invertible in $\cL(\cX)$ for all $\Delta \in
       \overline \cB \boldsymbol \Delta : = \{ \Delta \in \boldsymbol \Delta
       \colon \| \Delta \| \le 1 \}$.

       \item The LFT model $\Sigma$ is {\em $\cQ$-stable} if there exists an
invertible $Q
       \in \cD_{\boldsymbol \Delta}$ so that $\| Q^{-1} A Q \| < 1$, or,
       equivalently, if
       there exists a (strictly) positive-definite $X \in
       \cD_{\boldsymbol \Delta}$ (written
       as $X > 0$) so that
       $A X A^{*} - X < 0$.

       \item The LFT model $\Sigma$ has {\em robust performance} if
       $\Sigma$ is robustly stable and if in addition
       $ \| \cF_{u}\left( \sbm{A & B \\ C & D}, \Delta \right) \| < 1$ for
       all $\Delta \in \overline \cB \boldsymbol \Delta$.

       \item The LFT model $\Sigma$ has {\em $\cQ$-performance} if there exists an
       invertible $Q \in \cD_{\boldsymbol \Delta}$ so that
       \begin{equation}  \label{Qstable1}
    \left\| \begin{bmatrix} Q^{-1} & 0 \\ 0 &
       I_{\cY} \end{bmatrix} \begin{bmatrix} A & B \\ C & D
   \end{bmatrix} \begin{bmatrix} Q & 0 \\ 0 &  I_{\cU}
\end{bmatrix} \right\| < 1,
\end{equation}
or, equivalently, if there exists an $X > 0$ in $\cD_{\boldsymbol
\Delta}$
   so that
\begin{equation} \label{Qstable2}
\begin{bmatrix} A & B \\ C & D \end{bmatrix} \begin{bmatrix} X & 0 \\
      0 & I_{\cU} \end{bmatrix} \begin{bmatrix} A & B \\ C & D
\end{bmatrix}^* - \begin{bmatrix} X & 0 \\ 0 &  I_{\cY}
\end{bmatrix} < 0.
\end{equation}
   \end{enumerate}

   \begin{remark}  \label{R:clarify}  {\em
       A couple of remarks are in order
       to clarify these definitions.
       \begin{enumerate}
      \item[(i)] Note that the {\em robust stability} condition (1) and
      the {\em $\cQ$-stability} condition (2) involve only the
      operator $A \colon \cX \to \cX$.  In particular, one can
      replace the system matrix $\sbm{ A & B \\ C & D }$ by
      $$
        \begin{bmatrix} A & 0 \\ 0 & 0 \end{bmatrix} \colon
            \begin{bmatrix} \cX \\ \{0\} \end{bmatrix} \to
            \begin{bmatrix} \cX \\ \{0\} \end{bmatrix}
     $$
     without affecting the robust stability or $\cQ$-stability of the
     LFT model.

\item[(ii)]
Robust performance implies robust stability by definition. It is less
obvious but
also the case that $\cQ$-performance implies $\cQ$-stability. Given
that $\Sigma$ has
$\cQ$-performance, thus given an invertible $Q \in
\cD_{\boldsymbol\Delta}$ satisfying
(\ref{Qstable1}), it follows in particular that the upper left-hand corner of
the matrix inside the norm sign in (\ref{Qstable1}) also has norm
strictly less than
1, i.e., $ \| Q^{-1} A Q \| < 1$, which implies $\cQ$-stability.
   \end{enumerate}
   }\end{remark}

  The following result is well known (see \cite{LZD, DP}).

   \begin{proposition}  \label{P:relations}
       Let $\Sigma: = \left( \sbm{ A & B \\ C & D}, \, \boldsymbol
\Delta \right)$
       be an LFT-model system.  Then the following implications concerning
       $\Sigma$ hold:
       \begin{enumerate}
      \item $\cQ$-stability $\Longrightarrow$ robust
      stability.

      \item $\cQ$-performance $\Longrightarrow$ robust
      performance.
   \end{enumerate}
   Moreover, neither of the implications (1) nor (2) is reversible in
    general.
   \end{proposition}

%

For convenience in the discussion to follow, we assume
the input space and output space to be of the same
finite dimension $N$, and, in fact, make the identification
\[
\cU=\cY.
\]
The results can be extended to the case $\dim \cU\not=\dim \cY$ by using
the more general formalism of \cite{BGM3}.
We now specify the {\em full structure} $\boldsymbol \Delta_{\text{full}}$ by
\[
\boldsymbol \Delta_{\text{full}} = \cL(\cU,\cY)=\cL(\cU).
\]
 Note that
   $$
   \cD_{\boldsymbol \Delta_{\text{full}}} = \{ \lambda I_{\cU} \colon
   \lambda \in {\mathbb F}\}.
   $$
   We shall have use of the structure $\boldsymbol \Delta \oplus
   \boldsymbol \Delta_{\text{full}} \subset \cL(\cX \oplus \cU,\cX\oplus\cY)$
   consisting of operators of the form
   $$
   \boldsymbol \Delta \oplus \boldsymbol \Delta_{\text{full}} = \left\{
   \begin{bmatrix} \Delta & 0 \\ 0 & \Delta_{0} \end{bmatrix} \in
       \cL(\cX \oplus \cU,\cX\oplus\cY): \Delta \in \boldsymbol \Delta,\,  \Delta_{0}
       \in \Delta_{\text{full}} \right\}
    $$
    with associated commutant $\cD_{\boldsymbol \Delta \oplus
    \boldsymbol \Delta_{\text{full}}}$ given by
    $$
    \cD_{\boldsymbol \Delta \oplus
      \boldsymbol \Delta_{\text{full}}} = \left\{ \begin{bmatrix} Q & 0 \\ 0
      & \lambda I_{\cU} \end{bmatrix} \colon Q \in \cD_{\boldsymbol
      \Delta}, \, \lambda \in {\mathbb F} \right\}.
    $$

Note that the LFT model $\Sigma$ has {\em $\cQ$-performance} if and only if
       the system matrix $\sbm{A & B \\ C & D }$ is $\cQ$-stable
       with respect to the structure $\boldsymbol \Delta \oplus
       \boldsymbol \Delta_{\text{full}}$. Indeed, the condition that
   $\sbm{A & B \\ C & D }$ be $\cQ$-stable with respect to $\boldsymbol
   \Delta \oplus \boldsymbol \Delta_{\text{full}}$ a priori means that
   there exist an invertible $Q \in \cD_{\boldsymbol \Delta}$ and a
   nonzero number $\lambda$ so that
   \begin{equation}  \label{Qstable1'} \left\| \begin{bmatrix} Q^{-1} & 0 \\ 0 &
        \lambda^{-1} I_{\cU} \end{bmatrix} \begin{bmatrix} A & B \\ C & D
    \end{bmatrix} \begin{bmatrix} Q & 0 \\ 0 &  \lambda I_{\cU}
   \end{bmatrix} \right\| < 1,
   \end{equation}
   or, equivalently, there exist $X > 0$ in $\cD_{\boldsymbol \Delta}$
    and a number $\mu > 0$ so that
   \begin{equation}  \label{Qstable2'}
   \begin{bmatrix} A & B \\ C & D \end{bmatrix} \begin{bmatrix} X & 0 \\
       0 & \mu I_{\cU} \end{bmatrix} \begin{bmatrix} A & B \\ C & D
   \end{bmatrix}^* - \begin{bmatrix} X & 0 \\ 0 & \mu  I_{\cU}
   \end{bmatrix} < 0.
   \end{equation}
   However we can always replace $\sbm{ Q & 0 \\ 0 & \lambda I_{\cU} }$  by
   $\sbm {Q' & 0 \\ 0 & I_{\cU} } = \sbm{ \lambda^{-1} Q & 0 \\ 0 &
   I_{\cU}}$ in \eqref{Qstable1'} and $\sbm{ X & 0 \\ 0 & \mu I_{\cU}}$
   by $\sbm{X' & 0 \\ 0 & I_{\cU}} = \sbm{\mu^{-1} X & 0 \\ 0 &  I_{\cU}}$
   in \eqref{Qstable2'} to arrive at conditions of the respective forms
   \eqref{Qstable1} and \eqref{Qstable2}.  We shall see more of these
   simplifications via scaling in the sequel.

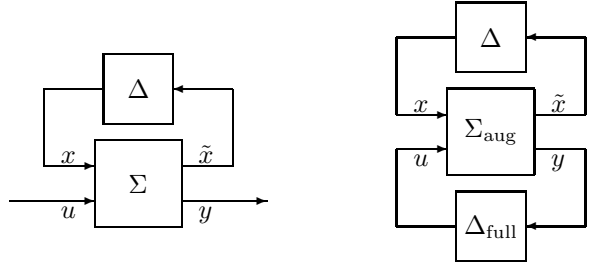
\begin{figure}[h]
\setlength{\unitlength}{0.09in}
\centering
\begin{picture}(10,14)
\put(-2,2){\framebox(5,5){$\Sigma$}}
\put(-1.5,8){\framebox(4,4){$\Delta$}}
\put(4,5.7){$\tilde x$}
\put(3,5.5){\line(1,0){3}}
\put(6,5.5){\line(0,1){4.5}}
\put(6,10){\vector(-1,0){3.5}}
\put(-4,5.7){$x$}
\put(-1.5,10){\line(-1,0){3.5}}
\put(-5,10){\line(0,-1){4.5}}
\put(-5,5.5){\vector(1,0){3}}
\put(4,2.5){$y$}
\put(3,3.5){\vector(1,0){5}}
\put(-4,2.5){$u$}
\put(-7,3.5){\vector(1,0){5}}
\end{picture}
\begin{picture}(10,14)
\put(8,5){\framebox(5,5){$\Sigma_{\text{aug}}$}}
\put(8.5,11){\framebox(4,4){$\Delta$}}
\put(8.5,0){\framebox(4,4){$\Delta_{\text{full}}$}}
\put(14,8.7){$\tilde x$}
\put(13,8.5){\line(1,0){3}}
\put(16,8.5){\line(0,1){4.5}}
\put(16,13){\vector(-1,0){3.5}}
\put(6,8.7){$x$}
\put(8.5,13){\line(-1,0){3.5}}
\put(5,13){\line(0,-1){4.5}}
\put(5,8.5){\vector(1,0){3}}
\put(14,5.5){$y$}
\put(13,6.5){\line(1,0){3}}
\put(16,6.5){\line(0,-1){4.5}}
\put(16,2){\vector(-1,0){3.5}}
\put(6,5.5){$u$}
\put(8.5,2){\line(-1,0){3.5}}
\put(5,2){\line(0,1){4.5}}
\put(5,6.5){\vector(1,0){3}}
\end{picture}
\caption{Original and augmented LFT-model}
\label{fig:augment}
\end{figure}

     Robust performance (or $\cQ$-performance) can be seen as simply
     robust stability (respectively, $\cQ$-stability) with respect to
     the appropriately contrived uncertainty structure (see Figure
    \ref{fig:augment}), as explained in the following proposition.

    \begin{proposition}  \label{P:stable=perform}
        Suppose that we are given an LFT model
        $$
        \Sigma = \left( \begin{bmatrix} A & B \\ C & D \end{bmatrix}
        \colon \begin{bmatrix} \cX \\ \cU \end{bmatrix} \to
        \begin{bmatrix} \cX \\ \cU \end{bmatrix}, \, \boldsymbol \Delta
       \right).
     $$
     Form the augmented LFT model $\Sigma_{\text{aug}}$ given by
     $$
     \Sigma_{\text{\rm aug}} = \left( \begin{bmatrix} \begin{bmatrix} A & B
     \\ C & D \end{bmatrix} & 0 \\ 0 & 0 \end{bmatrix} \colon
     \begin{bmatrix} \begin{bmatrix} \cX \\ \cU \end{bmatrix} \\ \{0\}
         \end{bmatrix} \to
         \begin{bmatrix} \begin{bmatrix}
       \cX \\ \cU \end{bmatrix} \\ \{0\}  \end{bmatrix},
       \boldsymbol \Delta \oplus \boldsymbol \Delta_{\text{\rm full}}
      \right).
     $$
     Then:
     \begin{enumerate}
         \item $\Sigma$ has robust performance with respect to
         $\boldsymbol \Delta$ if and only if
         $\Sigma_{\text{\rm aug}}$ is robustly stable with respect to
         $\boldsymbol \Delta \oplus \boldsymbol \Delta_{\text{\rm full}}$.

         \item $\Sigma$ has $\cQ$-performance with respect to
         $\boldsymbol \Delta$
         if and only if $\Sigma_{\text{\rm aug}}$ is $\cQ$-stable
         with respect to $\boldsymbol \Delta \oplus \boldsymbol
         \Delta_{\text{\rm full}}$.
         \end{enumerate}
     \end{proposition}

     \begin{proof}  See \cite{DWS, LZD}.
        \end{proof}

   The general philosophy of feedback control is: {\em
   given a plant with deficient properties} (e.g., lack of stability
   or performance), {\em design a compensator so that these deficiencies
   are rectified in the resulting closed-loop system}.
 To this end, we
   suppose that we are given an LFT model with input space $\cU$ and
   output space $\cY$ having direct-sum decompositions
   $$ \cU = \begin{bmatrix} \cU_{1} \\ \cU_{2} \end{bmatrix}, \qquad
     \cY = \begin{bmatrix}  \cY_{1} \\ \cY_{2} \end{bmatrix}.
   $$
   Usually the spaces $\cU_{1}$, $\cU_{2}$, $\cY_{1}$ and $\cY_{2}$ have
   physical interpretations as disturbance, control, error and
   measurement signals respectively.
 Then the LFT model $\Sigma$ has
   the more detailed form
   \begin{equation}  \label{Sigma-detailed}
   \Sigma = \left( \left[ \begin{array}{c|cc} A & B_{1} & B_{2} \\
   \hline C_{1} & D_{11} & D_{12} \\ C_{2} & D_{21} & D_{22}
\end{array} \right] \colon
   \begin{bmatrix} \cX \\ \hline  \cU_{1} \\ \cU_{2} \end{bmatrix} \to
       \begin{bmatrix} \cX \\ \hline \cY_{1} \\ \cY_{2} \end{bmatrix},
      \boldsymbol \Delta \right).
    \end{equation}

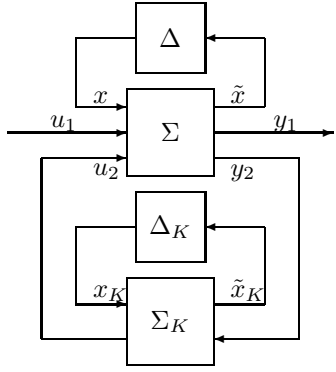
\begin{figure}[h]
\setlength{\unitlength}{0.09in}
\centering
\begin{picture}(10,21)
\put(2,11){\framebox(5,5){$\Sigma$}}
\put(2.5,17){\framebox(4,4){$\Delta$}}
\put(2,0){\framebox(5,5){$\Sigma_K$}}
\put(2.5,6){\framebox(4,4){$\Delta_K$}}
\put(8,15.2){$\tilde x$}
\put(7,15){\line(1,0){3}}
\put(10,15){\line(0,1){4}}
\put(10,19){\vector(-1,0){3.5}}
\put(0,15.2){$x$}
\put(2.5,19){\line(-1,0){3.5}}
\put(-1,19){\line(0,-1){4}}
\put(-1,15){\vector(1,0){3}}
\put(-2.5,13.8){$u_1$}
\put(-5,13.5){\vector(1,0){7}}
\put(10.5,13.8){$y_1$}
\put(7,13.5){\vector(1,0){7}}
\put(8,11){$y_2$}
\put(7,12){\line(1,0){5}}
\put(12,12){\line(0,-1){10.5}}
\put(12,1.5){\vector(-1,0){5}}
\put(0,11){$u_2$}
\put(2,1.5){\line(-1,0){5}}
\put(-3,1.5){\line(0,1){10.5}}
\put(-3,12){\vector(1,0){5}}
\put(8,3.9){$\tilde x_K$}
\put(7,3.5){\line(1,0){3}}
\put(10,3.5){\line(0,1){4.5}}
\put(10,8){\vector(-1,0){3.5}}
\put(0,3.9){$x_K$}
\put(2.5,8){\line(-1,0){3.5}}
\put(-1,8){\line(0,-1){4.5}}
\put(-1,3.5){\vector(1,0){3}}
\end{picture}
\caption{Feedback}
\label{fig:feedback}
\end{figure}

   Let us suppose that $\Sigma_{K}$
   is another LFT model of the form
   \begin{equation}  \label{Sigma-K}
   \Sigma_{K} = \left( \begin{bmatrix} A_{K} & B_{K} \\ C_{K} & D_{K}
   \end{bmatrix} \colon \begin{bmatrix} \cX_{K} \\ \cY_{2} \end{bmatrix}
   \to \begin{bmatrix} \cX_{K} \\ \cU_{2} \end{bmatrix}, \boldsymbol
   \Delta_{K} \right).
   \end{equation}
   Here the uncertainty structure $\boldsymbol \Delta_{K}$ for
   $\Sigma_{K}$ may be independent of the uncertainty structure
   $\boldsymbol \Delta$ for the original LFT model $\Sigma$ but we will
   be primarily interested in the case where there is a coupling
   between $\boldsymbol \Delta$ and $\boldsymbol
   \Delta_{K}$: we shall give a concrete model for this setup below.
   In any case, we may form the feedback connection
   \begin{align*}
     &  \begin{bmatrix} A & B_{1} & B_{2} \\ C_{1} & D_{11} & D_{12} \\
       C_{2} & D_{21} & D_{22} \end{bmatrix} \begin{bmatrix} x \\ u_{1}
       \\ u_{2} \end{bmatrix} = \begin{bmatrix} \widetilde x \\ y_{1} \\ y_{2}
   \end{bmatrix}, \qquad
   \begin{bmatrix} A_{K} & B_{K} \\ C_{K} & D_{K} \end{bmatrix}
       \begin{bmatrix} x_{K} \\ y_{2} \end{bmatrix} = \begin{bmatrix}
      \widetilde x_{K} \\ u_{2} \end{bmatrix}
    \end{align*}
    with transfer function $(\Delta, \Delta_{K}) \mapsto G_{cl}(\Delta,
    \Delta_{K})$ obtained by imposing the additional feedback equations
    $ x = \Delta \widetilde x$, $x_{K} = \Delta_{K} \widetilde x_{K}$
    (see Figure \ref{fig:feedback}).

    The resulting closed-loop transfer function is then given by
    $$
    G_{cl}(\Delta, \Delta_{K}) =
    \cF_{\ell}\left( \begin{bmatrix} G_{11}(\Delta) &
    G_{12}(\Delta) \\
    G_{21}(\Delta) & G_{22}(\Delta) \end{bmatrix},
    \cF_{u}\left(
    \begin{bmatrix} A_{K} & B_{K} \\ C_{K} & D_{K} \end{bmatrix}, \,
        \Delta_{K} \right) \right)
        $$
     where
    $$\begin{bmatrix} G_{11}(\Delta) & G_{12}(\Delta) \\
      G_{21}(\Delta) & G_{22}(\Delta) \end{bmatrix}
      = \cF_{u}\left(
      \left[ \begin{array}{c|cc} A & B_{1} & B_{2} \\ \hline
      C_{1} & D_{11} & D_{12} \\ C_{2} & D_{21} & D_{22} \end{array}
      \right],  \, \Delta  \right).
$$
As has been observed in \cite{Scherer, IS} and elsewhere
(at least for the case where $D_{22} = 0$), one can
realize $G_{cl}(\Delta, \Delta_{K})$ directly as the transfer
function of a linear-fractional model
\begin{equation}  \label{cltransfunc}
    G_{cl}(\Delta, \Delta_{K}) = \cF_{u}\left( \begin{bmatrix} A_{cl} &
    B_{cl} \\ C_{cl} & D_{cl} \end{bmatrix}, \begin{bmatrix} \Delta & 0
    \\ 0 & \Delta_{K} \end{bmatrix} \right)
\end{equation}
where the closed-loop state matrix $\sbm{A_{cl} & B_{cl} \\ C_{cl} &
D_{cl} }$ is given by
\begin{equation}  \label{clsys}
\begin{bmatrix} A_{cl} & B_{cl} \\ C_{cl} & D_{cl} \end{bmatrix} =
\cF_{\ell}\left( \left[ \begin{array}{ccc|cc}
A & 0 & B_{1} & 0 & B_{2}  \\
0 & 0 & 0 &  I & 0 \\
C_{1} & 0 & D_{11} &  0 & D_{12}  \\
\hline 0 & I & 0 & 0 & 0  \\
C_{2} & 0 & D_{21}  & 0 & D_{22} \\
\end{array} \right], \begin{bmatrix} A_{K} & B_{K} \\ C_{K} & D_{K}
\end{bmatrix} \right).
\end{equation}
The feedback-loop is well-posed exactly when $I - D_{22} D_{K}$ is
invertible. Since, under the assumption of well-posedness, one can
always arrange,
via a change of variable on the input-output space, that $D_{22}=0$,
it is usually
assumed that $D_{22} = 0$; in this case well-posedness is automatic
and $\sbm{A_{cl} &
B_{cl} \\ C_{cl} & D_{cl}}$ can be written out explicitly as
\begin{equation}  \label{clsysmat}
   \begin{bmatrix}  A_{cl} & B_{cl} \\ C_{cl} & D_{cl} \end{bmatrix} =
       \left[ \begin{array}{cc|c}
       A+B_{2} D_{K}C_{2} & B_{2} C_{K} & B_{1} + B_{2} D_{K} D_{21} \\
       B_{K} C_{2} & A_{K} & B_{K} D_{21} \\
       \hline C_{1} + D_{12} D_{K} C_{2} & D_{12} C_{K} & D_{11} +
       D_{12} D_{K} D_{21} \end{array} \right].
\end{equation}

For the sequel it is convenient to assume $\dim\cU_1=\dim \cY_1$ and identify
$$
   \cU_{1} = \cY_{1}.
$$
We set $\boldsymbol \Delta_{\text{full}}$ equal to the {\em full
structure} on $\cU_{1} = \cY_{1}$, i.e.,
$$
   \boldsymbol \Delta_{\text{full}} =  \cL(\cU_{1}).
$$
Given such a pair of LFT models $\Sigma$ and $\Sigma_{K}$ as in
\eqref{Sigma-detailed} and \eqref{Sigma-K}, once we specify a
closed-loop structure ${\boldsymbol \Delta}_{cl}$
we make the following definitions:
\begin{enumerate}
      \item  The LFT-feedback system $(\Sigma, \Sigma_{K})$ is {\em
      robustly stable} if the closed-loop state matrix $A_{cl}$ is
      robustly stable with respect to $\boldsymbol \Delta_{cl}$:
    $$ I -  \Delta_{cl}  A_{cl} \text{ is invertible
   for all } \Delta_{cl} \in \overline \cB \boldsymbol \Delta_{cl}.
   $$

   \item The LFT-feedback system $(\Sigma, \Sigma_{K})$ is {\em
   $\cQ$-stable} if there exists an invertible $Q_{cl}  \in
   \cD_{\boldsymbol \Delta_{cl}}$  so that
   $\| Q_{cl}^{-1} A_{cl} Q_{cl} \| < 1$, or, equivalently, if there
   exists $X_{cl} \in \cD_{\boldsymbol \Delta_{cl}}$ so that
   $A_{cl} X_{cl} A_{cl}^{*} - X_{cl} < 0$.

   \item The LFT-feedback system $(\Sigma, \Sigma_{K})$ has {\em robust
   performance} if $(\Sigma, \Sigma_{K})$ is robustly stable and if in
   addition the closed-loop transfer function $G_{cl}$ given by
   \eqref{cltransfunc} satisfies
   \begin{equation}  \label{clperf}
       \| G_{cl}(\Delta_{cl}) \| < 1 \text{ for all }
       \Delta_{cl} \in
       \overline \cB(\boldsymbol \Delta_{cl}).
   \end{equation}

   \item The LFT-feedback system $(\Sigma, \Sigma_{K})$  has {\em
   $\cQ$-performance} if the closed-loop system matrix
   $\sbm{ A_{cl} & B_{cl} \\ C_{cl} & D_{cl} }$ is $\cQ$-stable with
   respect to the structure $\boldsymbol \Delta_{cl} \oplus \boldsymbol
   \Delta_{\text{full}}$, i.e., if there exists an invertible
   $ Q_{cl} \in \cD_{\boldsymbol \Delta_{cl}}$ so that
   $$
    \left\| \begin{bmatrix} Q_{cl}^{-1} & 0 \\ 0 & I_{\cU_{1}} \end{bmatrix}
       \begin{bmatrix} A_{cl} & B_{cl} \\ C_{cl} & D_{cl} \end{bmatrix}
      \begin{bmatrix} Q_{cl} & 0 \\ 0 & I_{\cU_{1}} \end{bmatrix}
          \right\| < 1,
    $$
    or, equivalently, if there exists $X_{cl} > 0$ in $\cD_{\boldsymbol
    \Delta_{cl}}$ so that
   $$
   \begin{bmatrix} A_{cl} & B_{cl} \\ C_{cl} & D_{cl} \end{bmatrix}
       \begin{bmatrix}  X_{cl} &
      0 \\ 0 &  I_{\cU_{1}} \end{bmatrix}
      \begin{bmatrix} A_{cl} & B_{cl} \\ C_{cl} & D_{cl}
      \end{bmatrix} ^{*} -
      \begin{bmatrix}  X_{cl} &  0 \\ 0 &  I_{\cU_{1}} \end{bmatrix}
          < 0.
   $$
   \end{enumerate}

   As a consequence of Proposition \ref{P:relations} and part (3) of Remark
   \ref{R:clarify} applied to the closed-loop system, we see that
   $\cQ$-stability for a feedback system $(\Sigma, \Sigma_{K})$ implies
   robust stability  and that $\cQ$-performance implies robust
   performance (even with $\cQ$-stability for the closed-loop system).
   We also note that the notion of $\cQ$-performance for a
   closed-loop system is equivalent to the controller $\Sigma_{K}$
solving the scaled
   $H^{\infty}$-problem as formulated in \cite{AG} (see Section
   \ref{S:LPV} below).

   Given an  LFT model $\Sigma$ of the form \eqref{Sigma-detailed}, the
   {\em robust stabilization problem} is to find an LFT feedback
   controller $\Sigma_K$ of the form \eqref{Sigma-K} so that the closed-loop
   system is robustly stable, while the {\em robust $H^{\infty}$-problem}
   is to find $\Sigma_{K}$ of the form \eqref{Sigma-K} so that the
   closed-loop system has robust
   performance. The $\cQ$-version of these problems is to find
   $\Sigma_{K}$ which achieves $\cQ$-stability and $\cQ$-performance,
   respectively, for the closed-loop system.  It happens
   that necessary and sufficient conditions for the existence of such a
   $\Sigma_{K}$ are only available in general for the $\cQ$-version of
   the problem (see \cite{LZD}, \cite{AG}); these same conditions then
   give sufficient conditions for the non-$\cQ$ versions of the problems.

    We assume that the controller also has the form of an LFT model.
    Thus
   controller state space $\cX_{K}$ has the form
   \begin{equation} \label{K}
   \cX_{K} = {\mathbb F}^{Z_{K}} \text{  and }
   \cX_{K} =  \oplus_{k=1}^{d} \cX_{K,k} \text{ with }
   \cX_{K,k} =  {\mathbb F}^{n_{Kk} \cdot m_{Kk}}
   \end{equation}
   with block structure of the form
   \begin{equation}  \label{K-structure}
    \boldsymbol \Delta_{K} = \left\{ \Delta =
    \begin{bmatrix} \Delta_{K1} & & \\ & \ddots & \\
        & & \Delta_{Kd}  \end{bmatrix} \colon
   \Delta_{K,k} = \Delta^{0}_{Kk} \otimes I_{n_{Kk}} \right\}.
   \end{equation}
   In addition we assume that the structure
   $\boldsymbol \Delta_{cl}$ for the closed-loop system involves a
   coupling between the structure for the open-loop plant and that of the
   controller given by $\Delta^{0}_{k} = \Delta^{0}_{Kk}$ for $k = 1,
   \dots, d$, i.e.,
   \begin{equation}  \label{clcoupled}
     \boldsymbol \Delta_{cl}  = \left\{
     \begin{bmatrix}
         \Delta^{0}_{1} \otimes I_{n_{1}} & & & & & \\
         & \ddots & & & & \\
         & & \Delta^{0}_{d} \otimes I_{n_{d}} & & & \\
         & & & \Delta^{0}_{1} \otimes I_{n_{K1}} & & \\
          & & & & \ddots & \\
     & & & & & \Delta^{0}_{d} \otimes I_{n_{Kd}} \end{bmatrix}
     \right\}
    \end{equation}
   where the same $m_{k} \times m_{k}$ matrix $\Delta^{0}_{k}$ appears
   in the plant block (with multiplicity $n_{k}$) and in the controller
   block (but with multiplicity $n_{Kk}$). This additional assumption
   puts no real restriction on the generality of the method, as one can
   deny the controller (or the LFT system) to have access to certain
   blocks in the uncertainty structure simply by setting $n_{Kk}$ (or $n_k$)
   equal to zero. If we introduce the
   permutation matrix $P$ which shuffles the coordinates of the
   closed-loop state space according to the rule
   $$
   P \colon \begin{bmatrix} x_{1} \\ \vdots \\ x_{d} \\ x_{K1} \\
   \vdots \\ x_{Kd} \end{bmatrix} \mapsto
   \begin{bmatrix} x_{1} \\ x_{K1} \\ \vdots \\ x_{d} \\ x_{Kd}
   \end{bmatrix},
   $$
   then the closed-loop structure $\boldsymbol \Delta_{cl}$
   in the new coordinates is given by $ P^{*} \boldsymbol \Delta_{cl}
   P$ and has the same form as \eqref{scalar-blocks} but with $n_{k} +
   n_{Kk}$ in place of $n_{k}$.  In this representation the associated
   commutant $\cD_{P^{*} \boldsymbol \Delta_{cl} P}$ therefore has the
   form \eqref{Dcl1}, \eqref{Dcl2} with the index $n_{k} + n_{Kk}$ in
   place of $n_{k}$.

   The definition of robust
   stability, $\cQ$-stability, robust performance and $\cQ$-perform\-ance
   we now take with respect to the coupled closed-loop
   structure given by \eqref{clcoupled}.  For the rest of the paper
   we assume that we are given a pair of LFT models $(\Sigma,
   \Sigma_{K})$ with this structure.

   With these preliminaries out of the way we can state the following
   precise result.

   \begin{theorem}  \label{T:3}  {\em ($\cQ$-performance via
      dynamic multidimensional output feedback: see Theorem 11.5
       in \cite{DP})}
       There exists a multidimensional dynamic feedback controller
       $\Sigma_{K}$
       $$
         \Sigma_{K} = \left( \sbm{ A_{K} & B_{K} \\ C_{K} & D_{K}},
        \boldsymbol \Delta_{K}\right)
        $$
    with coupled uncertainty structure as in \eqref{K} and
    \eqref{K-structure} and
    prescribed
       controller dimension indices $n_{1}$, $\dots$, $n_{K}$
       so that the closed-loop
       system $(\Sigma, \Sigma_{K})$ (with closed-loop block
       structure as in \eqref{clcoupled}) has $\cQ$-performance if
       and only if there exist positive-definite matrices $X, Y \in
       \cD_{\boldsymbol \Delta}$ so that
       \begin{align}
     & \begin{bmatrix} N_{c} & 0 \\ 0 & I \end{bmatrix} ^{*}
          \begin{bmatrix} AYA^{*} - Y & AYC_{1}^{*} & B_{1} \\
       C_{1}Y A^{*} & C_{1} Y C_{1}^{*} - I & D_{11} \\
     B_{1}^{*} & D_{11}^{*} & -I \end{bmatrix} \label{Y-LMI}
       \begin{bmatrix} N_{c} & 0 \\ 0 & I \end{bmatrix} <  0, \\
       & \begin{bmatrix} N_{o} & 0 \\ 0 & I \end{bmatrix} ^{*}
      \begin{bmatrix} A^{*} X A - X & A^{*} X B_{1} & C_{1} ^{*} \\
     B_{1}^{*} X A & B_{1}^{*} X B_{1} - I & D_{11}^{*} \\
     C_{1} & D_{11} & -I \end{bmatrix}
     \begin{bmatrix} N_{o} & 0 \\ 0 & I \end{bmatrix} < 0,
         \label{X-LMI}
    \end{align}
      where $N_{c}$ and $N_{o}$ are matrices chosen so that
      \begin{align*}
      & N_{c} \text{ is injective and }
       \operatorname{Im } N_{c} = \operatorname{Ker } \begin{bmatrix}
       B_{2}^{*} & D_{12}^{*} \end{bmatrix} \text{ and } \\
       & N_{o} \text{ is injective and } \operatorname{Im } N_{o} =
       \operatorname{Ker } \begin{bmatrix} C_{2} & D_{21} \end{bmatrix},
       \end{align*}
      and, if we write
     $$X = \begin{bmatrix} X_{1} & & \\ & \ddots & \\ & &
      X_{d}\end{bmatrix}, \qquad Y = \begin{bmatrix} Y_{1} & & \\ &
      \ddots & \\ & & Y_{d} \end{bmatrix}
      $$
      with
      $$
      X_{k} = \begin{bmatrix} X_{k,0} & & \\ & \ddots  & \\ & &
      X_{k,0} \end{bmatrix}, \qquad
      Y_{k} = \begin{bmatrix} Y_{k,0} & & \\ & \ddots  \\ & & Y_{k,0}
      \end{bmatrix}
      $$
       as in the representation \eqref{Dcl1} and \eqref{Dcl2} for
       $\cD_{\boldsymbol \Delta}$,
       then we also have
      \begin{equation}  \label{coupling}
      \begin{bmatrix} X_{k,0} & I \\ I & Y_{k,0} \end{bmatrix} \ge 0
     \text{ and } \operatorname{rank } \begin{bmatrix} X_{k,0} & I
     \\ I & Y_{k,0} \end{bmatrix} \le n_{k} + n_{Kk}
      \text{ for }  k  = 1, \dots, d.
      \end{equation}
  \end{theorem}

     A special case of Theorem \ref{T:3} is the case where one insists
     that the controller be static, i.e., that all the controller
     state-space dimensions $n_{K1}, \dots, n_{Kd}$ be equal to $0$.  In
     this case, via a Schur-complement argument, one can see that
     the coupling condition \eqref{coupling} assumes the simple form
     \begin{equation}  \label{coupling'}
         Y_{k,0} = X_{k,0}^{-1} \text{ for } k = 1,
\dots,d,\quad\text{i.e.,}\quad Y=X^{-1}.
      \end{equation}
      We remark that the paper \cite{AG} as well as the exposition in
      the book \cite{DP} arrive at Theorem \ref{T:3} directly while the
      paper \cite{IS} (see also \cite{SIG}),
      explicitly only for the case $d=1$ but with an argument
      extendable to the general case here,
      first prove the special case for a static controller
      (conditions \eqref{Y-LMI}, \eqref{X-LMI} and \eqref{coupling'})
      and then use the observation \eqref{clsys} to reduce the
      dynamic-controller case to the static-controller case.

     \section{$\cQ$-stabilization as a consequence of closed-loop
$\cQ$-performance
     via feedback}\label{sec:Qstable}

     By zeroing out the disturbance and error channels, any
     $\cQ$-performance result leads to a $\cQ$-stability result.
     Application of this simple idea to Theorem \ref{T:3} leads to
     the following $\cQ$-stabilization result which we have not
     seen stated explicitly in the literature.

     \begin{theorem} \label{T:4} {\em ($\cQ$-stabilization via
     dynamic multidimensional output feedback: prescribed
     controller state-space dimensions)}
     There exists a multidimensional dynamic output controller
      $\Sigma_{K} = \left( \sbm{ A_{K} & B_{K} \\ C_{K} & D_{K}},
      \boldsymbol \Delta_{K}\right)$ (i.e., as in \eqref{K} and
      \eqref{K-structure} with prescribed
      dimension indices $n_{K1}, \dots, n_{Kd}$) so that the closed-loop
      system $(\Sigma, \Sigma_{K})$ (with closed-loop block structure
      $\boldsymbol \Delta_{cl}$ as in \eqref{clcoupled})
      is $\cQ$-stable if and only if
      there exist positive-definite matrices $X \in \cD_{\boldsymbol
      \Delta}$ and $Y \in \cD_{\boldsymbol \Delta}$ which satisfy the
      following pair of LMIs:
  \begin{align}
        & B_{\perp}^{*} A Y A^{*} B_{\perp} - B_{\perp}^{*} Y
         B_{\perp} < 0, \label{Y-LMI:T4} \\
       & C_{\perp} A^{*} X A C_{\perp}^{*} - C_{\perp} X
       C_{\perp}^{*} < 0,
       \label{X-LMI:T4}
     \end{align}
     where the matrices $B_{\perp}$ and $C_{\perp}$  are chosen so that
     \begin{align*}
        & B_{\perp} \text{ is injective and }
         \operatorname{Im } B_{\perp} = \operatorname{Ker } B_{2}^{*}, \\
       & C_{\perp}^{*} \text{ is surjective  and }  \operatorname{Im }
       C_{\perp}^{*} = \operatorname{Ker } C_{2}.
     \end{align*}
      Here $X$ and $Y$ have the block diagonal form as in \eqref{Dcl1}
      and \eqref{Dcl2}
      \begin{align}
      X = \begin{bmatrix} X_{1} & & \\ & \ddots & \\ & & X_{d}
       \end{bmatrix} &\text{ with } X_{k} = \begin{bmatrix} X_{k,0} & &
       \\ & \ddots & \\ & & X_{k,0} \end{bmatrix}, \label{Xdiag}\\
        Y = \begin{bmatrix} Y_{1} & & \\ & \ddots & \\ & & Y_{d}
    \end{bmatrix} &\text{ with } Y_{k} = \begin{bmatrix} Y_{k,0} & &
    \\ & \ddots & \\ & & Y_{k,0} \end{bmatrix} \label{Ydiag}
    \end{align}
      and must in addition satisfy the coupling and rank conditions
      \begin{equation}  \label{ranks}
      \begin{bmatrix} X_{k,0} & I \\ I & Y_{k,0} \end{bmatrix}
          \ge 0 \text{ and }
          \operatorname{rank } \begin{bmatrix} X_{k,0} & I \\ I &
          Y_{k,0} \end{bmatrix} \le n_{k} + n_{Kk} \text{ for } k =
          1, \dots, d.
       \end{equation}
       \end{theorem}

       \begin{proof}  It suffices to apply the observation (i) in Remark
           \ref{R:clarify} to the closed-loop system and set the input
           space $\cU_{1}$ and output space $\cY_1$ equal to
           $\{0\}$ in Theorem \ref{T:3}. Note that in this case the matrices
           $N_c$ and $N_o$ in Theorem \ref{T:3} coincide with $B_\perp$
           and $C_\perp^*$, respectively.
       \end{proof}

There are two extreme special cases of Theorem \ref{T:4}: (1) the
case where we prescribe $n_{Kk}=0$ for each $k = 1, \dots, d$, and
(2) the case where no bounds are imposed on $n_{Kk}$.  In each of
these cases, the coupling and rank conditions \eqref{ranks} either
disappear or can be put in a different form.  In this way we recover
$\cQ$-stabilization results appearing in \cite{LZD} as special cases.

\begin{theorem}  \label{T:1and2} \textbf{
      {\em (1) $\cQ$-stabilization via static output
      feedback (see Theorem III-9 in \cite{LZD}:}})
      There exists a static output feedback controller (i.e.,
      $\Sigma_{K} = (D_{K}, 0)$ where $n_{Kk} = 0 $ for $k = 1, \dots,
      d$) so that the closed-loop system $(\Sigma, \Sigma_{K})$ is
      $\cQ$-stable if and only if there exists a positive-definite
      matrix $X \in \cD_{\boldsymbol \Delta}$ so that the following
      two LMIs hold:
    \begin{align}
       & B_{\perp}^{*} A X^{-1} A^{*} B_{\perp} - B_{\perp}^{*} X^{-1}
    B_{\perp} < 0, \label{Y-LMI:T1} \\
      & C_{\perp} A^{*} X A C_{\perp}^{*} - C_{\perp} X
      C_{\perp}^{*} < 0.
      \label{X-LMI:T1}
    \end{align}
    Here the matrices $B_{\perp}$ and $C_{\perp}$  are chosen as in Theorem \ref{T:4}.

\medskip
\noindent
\textbf{\em (2) $\cQ$-stabilization via
      dynamic  multidimensional output feedback (see Theorem V-1 in
       \cite{LZD}):}
    There exists a multidimensional dynamic feedback controller
    $$
    \Sigma_{K} = \left( \sbm{ A_{K} & B_{K} \\ C_{K} & D_{K}},
       \boldsymbol \Delta_{K}\right)
       $$
       (i.e., as in \eqref{K} and
       \eqref{K-structure} with no restriction on the
       dimension indices $n_{K1}$, $\dots$, $n_{Kd}$) so that the closed-loop
       system $(\Sigma, \Sigma_{K})$ (with closed-loop block structure
       $\boldsymbol \Delta_{cl}$ as in \eqref{clcoupled})
       is $\cQ$-stable if and only if
       there exist positive-definite matrices $X \in \cD_{\boldsymbol
       \Delta}$ and $Y \in \cD_{\boldsymbol \Delta}$ which satisfy the
       following pair of LMIs:
       \begin{align}
      & A Y A^{*} - Y - B_{2} B_{2}^{*} < 0, \label{Y-LMI:T2} \\
      & A^{*} X A - X - C_{2}^{*} C_{2} < 0.  \label{X-LMI:T2}
   \end{align}

    \end{theorem}

    \begin{proof}  To prove the first statement,
    apply Theorem \ref{T:4} to the case
      where $n_{K1} = \cdots = n_{Kd} = 0$.
      Note that the conditions \eqref{Y-LMI:T4} and \eqref{X-LMI:T4} are
      exactly conditions \eqref{Y-LMI:T1} and \eqref{X-LMI:T1} but
      with $Y$ taken to be equal to $X^{-1}$.
      Note also that the rank condition
      $$ \operatorname{rank}\, \begin{bmatrix} X_{k,0} & I \\ I & Y_{k,0}
       \end{bmatrix} = n_{k}
      $$
      (where $n_{k} = \operatorname{rank}\, X_{k,0} =
      \operatorname{rank}\, Y_{k,0}$ for each $k$) is equivalent to
      $X_{k,0}=Y_{k,0}^{-1}$ for each $k$, or $X = Y^{-1}$.

      To prove the second statement,
 apply Theorem \ref{T:4} to the case where there are
no restrictions on the dimension indices  $n_{K1}, \dots, n_{Kd}$.
Let $B_\perp$ and $C_\perp$ be as in Theorem \ref{T:4}.  As pointed
out by one of the reviewers, the existence of positive-definite
$X,Y \in \cD_{\boldsymbol \Delta}$ satisfying \eqref{Y-LMI:T2} and
\eqref{X-LMI:T2} is equivalent to existence of (not necessarily the same)
positive-definite $X,Y \in \cD_{\boldsymbol \Delta}$ satisfying
\eqref{X-LMI:T4} and \eqref{Y-LMI:T4}, by a simple application of
Finsler's lemma (see \cite[Lemma 3]{IS}).

 As we are imposing no constraints on
the control state-space dimension indices $n_{K1}, \dots,
n_{Kd}$, the rank conditions in \eqref{ranks} can safely be
ignored.  To handle the coupling conditions
\begin{equation}  \label{coupling''}
\begin{bmatrix} X_{k,0} & I \\ I & Y_{k,0} \end{bmatrix} \ge 0,
\end{equation}
note that we can always replace $X>0$ and $Y>0$ by $\widetilde X = \mu X$, $\widetilde Y = \mu Y$
with the scalar multiplier $\mu > 0$ sufficiently large to guarantee
\eqref{coupling''} (with $\widetilde Y, \widetilde X$ in place of
$Y,X$)  while not affecting the validity of the homogeneous LMIs
\eqref{Y-LMI:T4} and \eqref{X-LMI:T4}.
\end{proof}

     \section{Closed-loop $\cQ$-performance as a consequence of
     $\cQ$-stabilization}\label{sec:Qperformance}

     In this section we give two illustrations of the
     Principle of Reduction of Robust Performance to Robust Stabilization
     given in the introduction.  Proposition \ref{P:stable=perform} is
     one such illustration, but note that Proposition
     \ref{P:stable=perform} pays no heed to compensator partial-state
     dimension.  Application of the idea in Proposition \ref{P:stable=perform}
     to Theorem 3.2 (2) leads to the following result.

     \begin{theorem}  \label{T:5} {\em ($\cQ$-performance via a
         dynamic output controller with input-out\-put-loop dynamics)}
        There exists a multidimensional dynamic output
        feedback-con\-troller as in the configuration on the left
    side of Figure \ref{fig:fullblock} which achieves
         $\cQ$-perfor\-mance for the closed-loop system if and
        only if there exist positive-definite matrices $X, Y \in
        \cD_{\boldsymbol \Delta}$ which satisfy the LMIs \eqref{Y-LMI}
        and \eqref{X-LMI}.

        \end{theorem}

    \begin{remark} \label{R:artifical}  {\em We emphasize that the
        feedback configuration on the left side of Figure
        \ref{fig:fullblock} is contrived and not of interest from
        the physical point of view.  The point here is that
        adherence to the Principle of Reduction of Robust Performance to Robust Stabilization
         does give the equivalence between two
        control problems, but sometimes not between problems of
        practical interest, contrary to expectations as suggested in
        \cite{LZD}.
         }\end{remark}

        \begin{proof}
    Let $\Sigma$ be the LFT model (\ref{Sigma-detailed}) and $\Sigma_K$ the
LFT model (\ref{Sigma-K}). By definition, the closed-loop LFT-feedback
system $\Sigma_{cl}=(\Sigma,\Sigma_K)$ has $\cQ$-performance if the
closed-loop system
matrix $\sbm{ A_{cl} & B_{cl} \\ C_{cl} & D_{cl} }$ in (\ref{clsysmat})
is $\cQ$-stable.

Next we introduce the adjusted LFT model $\Sigma_{adj}$ given by
\begin{equation}  \label{Sigma-adj}
   \Sigma_{adj} = \left( \left[
\begin{array}{cc|cc} A & B_{1} & 0 & B_{2} \\ C_{1} & D_{11} & 0 &
      D_{12} \\ \hline 0 & 0 & 0 & 0 \\ C_{2} & D_{21} & 0 & 0
\end{array} \right] \colon \begin{bmatrix} \begin{bmatrix} \cX \\ \cU_{1}
\end{bmatrix} \\ \{0\} \\ \cU_{2} \end{bmatrix} \to
\begin{bmatrix} \begin{bmatrix} \cX \\ \cU_{1}
\end{bmatrix} \\ \{0\} \\ \cU_{2} \end{bmatrix}, \boldsymbol \Delta
\oplus \boldsymbol \Delta_{\text{full}} \right),
\end{equation}
and its closed-loop LFT model $\Sigma_{adj,cl}=(\Sigma_{adj},\Sigma_K)$.

We claim that $\Sigma_{cl}$ has $\cQ$-performance if and only if
$\Sigma_{adj,cl}$
is $\cQ$-stable. To see this, note that the state operator
$A_{adj,cl}$ for the LFT model $\Sigma_{adj,cl}$ is given by
\begin{eqnarray*}
A_{adj,cl}&=&\mat{cc}
{\mat{cc}{A&B_1\\C_1&D_{11}}+\mat{c}{B_2\\D_{12}}D_K\mat{cc}{C_2&D_{2,1}}
&\mat{c}{B_2\\D_{12}}C_K\\B_K\mat{cc}{C_2&D_{21}}&A_K}\\
&=&\mat{ccc}{
A+B_2D_KC_2&B_1+B_2D_KD_{21}&B_2C_K\\
C_1+D_{12}D_KC_2&D_{11}+D_{12}D_KD_{21}&D_{12}C_K\\
B_KC_2&B_KD_{21}&A_K}.
\end{eqnarray*}
By rearranging rows and columns we can identify $A_{adj,cl}$ with the
closed-loop system
matrix $\sbm{ A_{cl} & B_{cl} \\ C_{cl} & D_{cl} }$ as in \eqref{clsysmat};
in particular, it follows that $A_{adj,cl}$ is $\cQ$-stable if and only if
$\sbm{ A_{cl} & B_{cl} \\ C_{cl} & D_{cl} }$ is $\cQ$-stable as claimed.

Applying Theorem \ref{T:1and2} (2) to the adjusted LFT model $\Sigma_{adj}$
thus provides us
with necessary and sufficient conditions for the existence of a
multidimensional feedback
controller $\Sigma_K$ so that the closed-loop system $\Sigma_{cl}$
has $\cQ$-performance
and that has access to $\boldsymbol \Delta_K$ in (\ref{K-structure})
as well as to the
full block $\boldsymbol \Delta_{\text{full}}$.

As was already remarked in the proof of Theorem \ref{T:1and2},
as a consequence of the Finsler lemma
the LMIs
(\ref{Y-LMI:T2}), (\ref{X-LMI:T2}) are equivalent to the LMIs
(\ref{Y-LMI:T4}),
(\ref{X-LMI:T4}). It thus remains to show that the LMIs
(\ref{Y-LMI:T4}) and (\ref{X-LMI:T4})
when specified to $\Sigma_{adj}$ are equivalent to the LMIs
(\ref{Y-LMI}) and (\ref{X-LMI}).
Notice that the matrices $B_\perp$ and $C_\perp$, when specified for
$\Sigma_{adj}$ rather than
for $\Sigma$, coincide with $N_c$ and $N_o^*$ in Theorem \ref{T:3}.
For the record we note that
\eqref{Y-LMI:T4} and \eqref{X-LMI:T4}, spelled out for the case at
hand, assume the form
\begin{align}
      N_{c}^{*} \left( \begin{bmatrix} A & B_{1} \\ C_{1} & D_{11}
\end{bmatrix} \begin{bmatrix} Y & 0 \\ 0 & \mu I \end{bmatrix}
\begin{bmatrix} A^{*} & C_{1} ^{*} \\ B_{1} ^{*} & D_{11} ^{*}
\end{bmatrix} - \begin{bmatrix} Y & 0 \\ 0 & \mu I \end{bmatrix} \right)
N_{c} < 0, \label{Y-LMI:T4''} \\
    N_{o}^{*} \left( \begin{bmatrix} A^{*} & C_{1}^{*} \\ B_{1}^{*} &
    D_{11}^{*} \end{bmatrix} \begin{bmatrix} X & 0 \\ 0 &
    \widetilde  \mu I \end{bmatrix}
    \begin{bmatrix} A & B_{1} \\ C_{1} & D_{11}
\end{bmatrix} - \begin{bmatrix} X & 0 \\ 0 & \widetilde \mu
I \end{bmatrix} \right)
        N_{o}< 0  \label{X-LMI:T4''}
\end{align}
with $N_{c}$ and $N_{o}$ as in Theorem \ref{T:3}.
As these inequalities are homogeneous in $\sbm{ Y & 0 \\ 0 & \mu I }$
and $\sbm{ X & 0 \\ 0 & \widetilde \mu I}$ respectively, at this
stage we may rescale if necessary to arrange without loss of
generality that $\mu = \widetilde \mu = 1$.  Theorem \ref{T:5}
follows once we see that conditions \eqref{Y-LMI:T4''} and
\eqref{X-LMI:T4''} can be converted to the more linear form of
conditions \eqref{Y-LMI} and \eqref{X-LMI}.

But this last step is a standard Schur-complement computation.  We
will show only that \eqref{Y-LMI} is equivalent to \eqref{Y-LMI:T4''}
as the equivalence of \eqref{X-LMI} with \eqref{X-LMI:T4''} is similar.
Rewrite \eqref{Y-LMI} in the form
$$
    \begin{bmatrix} N_{c}^{*} \begin{bmatrix} AYA^{*} - Y & A Y
      C_{1}^{*} \\ C_{1}Y A^{*} & C_{1} Y C_{1}^{*} - I \end{bmatrix}
      N_{c} & N_{c}^{*} \begin{bmatrix} B_{1} \\ D_{11} \end{bmatrix}  \\
      \begin{bmatrix} B_{1}^{*} & D_{11}^{*} \end{bmatrix} N_{c} & -I
       \end{bmatrix} < 0.
$$
Validity of \eqref{Y-LMI} is equivalent to
negative definiteness of the Schur complement with respect to the
lower right entry $-I$:
\begin{align*} 0 & > N_{c}^{*} \begin{bmatrix} AYA^{*} - Y & AYC_{1}^{*} \\
C_{1}Y A^{*} & C_{1} Y C_{1}^{*} -I \end{bmatrix} N_{c}
+ N_{c}^{*} \begin{bmatrix} B_{1} \\ D_{11} \end{bmatrix}
\begin{bmatrix} B_{1}^{*} & D_{11}^{*} \end{bmatrix} N_{c} \\
      & = N_{c}^{*} \begin{bmatrix} AYA^{*} - Y + B_{1} B_{1}^{*} &
      AYC_{1}^{*} + B_{1} D_{11}^{*} \\ C_{1}Y A^{*} + D_{11} B_{1}^{*}
      & C_{1} Y C_{1}^{*} - I + D_{11} D_{11}^{*} \end{bmatrix} N_{c}
      \end{align*}
   which, upon rearrangement, agrees with \eqref{Y-LMI:T4''} (with $\mu
   = 1$) as expected. This completes the proof of Theorem \ref{T:5}.
   \end{proof}

   We now show how imposing the condition that the controller
   state-space dimension constraint $n_{K,\cU_{1}}=0$ (see the right
   signal-flow diagram in Figure \ref{fig:fullblock}) leads to a proof
   of Theorem \ref{T:3} on $\cQ$-performance as a consequence of
   the result of Theorem \ref{T:4} on $\cQ$-stability;  there
   follows a presumably new interpretation of the coupling condition
   \eqref{coupling} in Theorem \ref{T:3} as the precise extra condition
   required in Theorem \ref{T:5} for the existence of a controller
$\Sigma_{K}$ as in
   Theorem \ref{T:5} which does not have access to the artificial
   full-block $\boldsymbol \Delta_{\text{full}}$.

        \begin{theorem} \label{T:6}  The $\cQ$-performance result
       Theorem \ref{T:3} can be seen as a
       corollary to the $\cQ$-stabilization result
       Theorem \ref{T:4}.
    \end{theorem}

    \begin{proof}  We follow the same scheme as used in the proof of
        Theorem \ref{T:5} above but now with use of Theorem \ref{T:4} rather
        than Theorem \ref{T:1and2} (2) and with the imposition of the
        constraint that the
        controller has no access to the artificial full block
        $\boldsymbol\Delta_{\text{full}}$. For the special situation where
        $\Sigma=\Sigma_{\text{adj}}$ as in \eqref{Sigma-adj},
        conditions \eqref{Y-LMI:T4} and \eqref{X-LMI:T4} in Theorem \ref{T:4}
        become the LMIs \eqref{Y-LMI:T4''} and \eqref{X-LMI:T4''}
        given above combined with the two coupling and rank conditions
        \begin{align}
      &    \begin{bmatrix} X_{k,0} & I \\ I & Y_{k,0} \end{bmatrix}
           \ge 0 \text{ and }
           \operatorname{rank}\,
           \begin{bmatrix} X_{k,0} & I \\ I & Y_{k,0} \end{bmatrix}
           \le n_{k} + n_{Kk}
            \text{ for } k = 1, \dots, d,
           \label{coupling1} \\
    & \begin{bmatrix} \mu I_{\cU_{1}} & I_{\cU_{1}} \\ I_{\cU_{1}} &
    \widetilde \mu I_{\cU_{1}} \end{bmatrix} \ge 0 \text{ and }
    \operatorname{rank} \,
    \begin{bmatrix} \mu I_{\cU_{1}} & I_{\cU_{1}} \\ I_{\cU_{1}} &
      \widetilde \mu I_{\cU_{1}} \end{bmatrix} = \operatorname{dim}\,
      \cU_{1},  \label{coupling2}
     \end{align}
     (where $X$ and $Y$ are given by the block diagonal forms as in
     \eqref{Xdiag} and \eqref{Ydiag}).

     As in the proof of Theorem \ref{T:5} we would like to rescale
     $\sbm{ Y & 0 \\ 0 & \mu I }$ and $\sbm{ X & 0 \\ 0 & \widetilde \mu I}$
     in \eqref{Y-LMI:T4''} and \eqref{X-LMI:T4''} with $\mu^{-1}$ and
     $\widetilde \mu^{-1}$, respectively, to obtain $\mu=\widetilde \mu=1$
     and arrive at the equivalence of conditions \eqref{Y-LMI:T4},
\eqref{X-LMI:T4}
     with conditions \eqref{Y-LMI}, \eqref{X-LMI}. First we need to
     check that
     this rescaling does not violate the coupling and rank conditions in
     \eqref{coupling1} and \eqref{coupling2}.

     The rank condition in \eqref{coupling2} forces
     $$
       \widetilde \mu = 1/\mu,
     $$
     and conversely, having $\widetilde \mu = 1/\mu$ the coupling and rank
     condition in \eqref{coupling2} are satisfied. Hence \eqref{coupling2}
     is equivalent to $\widetilde \mu = 1/\mu$; we therefore assume
     for the remainder of the proof that $\widetilde \mu = 1/\mu$.

     In particular, we may rescale
     $\sbm{ Y & 0 \\ 0 & \mu I }$ and $\sbm{ X & 0 \\ 0 & \widetilde \mu I}$
     without violating \eqref{coupling2}, but not independently: when we rescale
     $\sbm{ Y & 0 \\ 0 & \mu I }$ with
     $\alpha>0$, then $\sbm{ X & 0 \\ 0 & \widetilde \mu I}$ should be rescaled
     with $\alpha^{-1}$. In particular, taking $\alpha = \widetilde
     \mu$, given that $\widetilde \mu = 1/\mu$, leads to the desired
     result; after a rescaling we may take $\mu =
     \widetilde \mu = 1$ and maintain the validity of
     \eqref{Y-LMI:T4''}, \eqref{X-LMI:T4''} and \eqref{coupling2}.

     It remains to see whether \eqref{coupling1} still holds under
this rescaling.
     To verify this, observe that
     \begin{eqnarray*}
     \mat{cc}{\widetilde \mu^{-1}X_{k,0}&I\\I&\widetilde \mu Y_{k,0}}
     &=&\mat{cc}{\widetilde \mu^{-1}I&0\\0&I}\mat{cc}{X_{k,0}&I\\I&Y_{k,0}}
     \mat{cc}{I&0\\0&\widetilde \mu I}\\
     &=&\widetilde \mu^{-1}\mat{cc}{I&0\\0&\widetilde \mu
I}\mat{cc}{X_{k,0}&I\\I&Y_{k,0}}
     \mat{cc}{I&0\\0&\widetilde \mu I}.
     \end{eqnarray*}
     We thus obtain that the conditions \eqref{Y-LMI:T4''},
     \eqref{X-LMI:T4''}, \eqref{coupling1} and \eqref{coupling2} are
     exactly equivalent to the conditions \eqref{Y-LMI}, \eqref{X-LMI}
     and \eqref{coupling} given in Theorem \ref{T:3}, and we arrive at the
     $\cQ$-performance result Theorem \ref{T:3} as a consequence
     of the $\cQ$-stabilization result Theorem \ref{T:4} as asserted.
     \end{proof}

     \section{Applications}\label{sec:applications}
     In this section we discuss how the abstract results on LFT model
     systems of the previous section apply to more concrete control
     settings.
     We discuss two particular applications: {\em robust control for
     systems with structured uncertainty} and {\em robust control for
     LPV systems}.

     \subsection{Systems with LFT models for structured uncertainty}

     We suppose that we are given a standard linear time-invariant
     input/state/output linear system model
     $$\Sigma \colon \left\{ \begin{array}{ccc}
      x(t+1) & = & A_{M}(\delta_{U}) x(t)+ B_{M1}(\delta_{U}) w(t) +
      B_{M2}(\delta_{U}) u(t) \\
      z(t) & = & C_{M1}(\delta_{U}) x(t) + D_{M11}(\delta_{U}) w(t) +
      D_{M12}(\delta_{U}) u(t)  \\
      y(t) & = & C_{M2}(\delta_{U}) x(t) + D_{M21}(\delta_{U}) w(t)
      + D_{M22}(\delta_{U}) u(t)
      \end{array}  \right.
      $$
      where the system matrix
      $$\begin{bmatrix}  A_{M}(\delta_{U}) & B_{M1}(\delta_{U}) &
      B_{M2}(\delta_{U}) \\ C_{M1}(\delta_{U})  & D_{M11}(\delta_{U}) &
      D_{M12}(\delta_{U}) \\ C_{M2}(\delta_{U}) & D_{M21}(\delta_{U}) &
      D_{M22}(\delta_{U}) \end{bmatrix}
      $$
      is not known exactly but depends on some
      {\em uncertainty} parameters $\delta_{U}  = (\delta_{1}, \dots,
      \delta_{d})$.  Here the quantities $\delta_{i}$ are viewed as
         uncertainties unknown to the controller.  The goal is to design
         a controller $\Sigma_{K}$ (independent of $\delta_U$)
         so that the closed-loop system has desirable
         properties for all admissible values of $\delta_U$, usually
         normalized to be $|\delta_{k}| \le 1$ for $k = 1, \dots, d$.

      The transfer function for the uncertainty
      parameter $\delta_{U}$ can be expressed as
      \begin{align}
       G(\delta_{U}) & = \begin{bmatrix} D_{M11}(\delta_{U}) &
      D_{M12}(\delta_{U}) \\ D_{M21}(\delta_{U}) &
      D_{M22}(\delta_{U})\end{bmatrix}  \notag \\
     &   \qquad
     +  \lambda \begin{bmatrix} C_{M1}(\delta_{U}) \\ C_{M2}(\delta_{U})
\end{bmatrix} ( \lambda I_{n} - A_{M}(\delta_{U}))^{-1} \begin{bmatrix}
B_{M1}(\delta_{U}) & B_{M2}(\delta_{U}) \end{bmatrix} \notag \\
    &  = \cF_{u} \left( \left[ \begin{array}{c|cc}
       A_{M}(\delta_{U}) & B_{M1}(\delta_{U}) & B_{M2}(\delta_{U}) \\
       \hline C_{M1}(\delta_{U})  & D_{M11}(\delta_{U}) &
      D_{M12}(\delta_{U}) \\ C_{M2}(\delta_{U}) & D_{M21}(\delta_{U}) &
      D_{M22}(\delta_{U}) \end{array} \right], \lambda I \right)
      \label{Utransfunc}
     \end{align}
     where we have introduced the aggregate variable
     $$  \delta = (\delta_{U}, \lambda) = (\delta_{1}, \dots,
     \delta_{d}, \lambda).
     $$

      It is not too much of a restriction to assume in
      addition that the functional dependence on $\delta_{U}$ is given by a
      linear fractional map (where the subscript $U$ suggests {\em
      uncertainty} and the subscript $S$ suggests {\em shift})
     $$ \begin{bmatrix}  A_{M}(\delta_{U}) & B_{M1}(\delta_{U}) &
     B_{M2}(\delta_{U})
\\ C_{M1}(\delta_{U})  & D_{M11}(\delta_{U}) &
     D_{M12}(\delta_{U}) \\ C_{M2}(\delta_{U}) & D_{M21}(\delta_{U}) &
     D_{M22}(\delta_{U}) \end{bmatrix}   =
     \cF_{u}\left( \left[ \begin{array}{c|ccc}
     A_{UU} & A_{US} & B_{U1} & B_{U2} \\
     \hline A_{SU} & A_{SS} & B_{S1} & B_{S2} \\
     C_{1U} & C_{1S} & D_{11} & D_{12} \\
     C_{2U} & C_{2S} & D_{21} & D_{22} \end{array} \right],
     \Delta_{U} \right)
    $$
    where we take the uncertainty structure matrix $\Delta_{U}$
    to have the form as
    in \eqref{scalar-blocks} with $m_{k} = 1$ for $k = 1, \dots, d$
    for simplicity:
    $$  \Delta_{U} = \begin{bmatrix} \delta_{1} I_{n_{1}} & & \\ &
    \ddots & \\ & & \delta_{d} I_{n_{d}} \end{bmatrix}.
    $$
    Finally, if we introduce the aggregate matrix
    \begin{equation} \label{aggsys}
      \begin{bmatrix} A & B_{1} & B_{2} \\ C_{1} & D_{11} & D_{12} \\
    C_{2} & D_{21} & D_{22} \end{bmatrix} =
     \left[ \begin{array}{cc|cc}
       A_{UU} & A_{US} & B_{U1} & B_{U2} \\
       A_{SU} & A_{SS} & B_{S1} & B_{S2} \\
      \hline  C_{1U} & C_{1S} & D_{11} & D_{12} \\
       C_{2U} & C_{2S} & D_{21} & D_{22} \end{array} \right],
    \end{equation}
    then the transfer function $G(\delta)$ \eqref{Utransfunc}
    can conveniently be written in LFT form as
    $$
      G(\delta_{U}) = \cF_{u}\left( \left[ \begin{array}{c|cc}
      A & B_{1} & B_{2} \\ \hline C_{1} & D_{11} & D_{12} \\ C_{2} &
      D_{21} & D_{22} \end{array} \right], \Delta \right)
    $$
    where we have now set $\Delta$ equal to the expanded block diagonal matrix
    $$
      \Delta = \begin{bmatrix} \Delta_{U} & 0 \\ 0 & \lambda I
     \end{bmatrix}\quad(\lambda\in{\mathbb F}).
    $$

    If we take $\sbm{ A & B_{1} & B_{2} \\ C_{1} & D_{11} & D_{12} \\
    C_{2} & D_{21} & D_{22}}$ as in \eqref{aggsys} and introduce the
    block structure
    \begin{equation}  \label{aggblock}
     \boldsymbol \Delta = \left\{ \begin{bmatrix} \delta_{1}I_{n_{1}} &
     & & \\ & \ddots & & \\ & & \delta_{d} I_{n_{d}} \\ & & & \lambda
     I_{n} \end{bmatrix} \colon \delta_{1}, \dots, \delta_{d}, \lambda
     \in {\mathbb F} \right\}
    \end{equation}
    we may consider $(\Sigma, \boldsymbol \Delta)$ as an LFT model of
    the form \eqref{Sigma-detailed}.  Without much loss of generality,
    we follow the common normalization and assume that  $D_{22} = 0$.

    The problem is to design an output-feedback controller $K \colon y
    \mapsto u$ so that the closed-loop
    system
    $$ \left\{ \begin{array}{ccc}
     x(t+1) & = & A_{cl}(\delta_{U}) x(t) + B_{cl}(\delta_{U})  w(t) \\
     z(t) & = & C_{cl}(\delta_{U}) x(t) + D_{cl}(\delta_{U}) w(t)
     \end{array}
     \right.
     $$
     is {\em robustly stable} (i.e., $A_{cl}(\delta_{U})$ has spectral radius
     less than $1$ for all $\delta_{U}$ such that $|\delta_{k}| \le 1 $
     for each $k=1, \dots, d$) and, that perhaps also solves the {\em robust
     performance} problem, i.e., in addition the closed-loop transfer function
     $G(\delta)$ satisfies
     $$
       \| G(\delta) \| < 1 \text{ for all } \delta = (\delta_{1},
       \dots, \delta_{d}, \lambda) \text{ with } |\delta_{k}|, |\lambda| \le 1.
     $$

     If we only allow for static controllers, then a necessary and
     sufficient condition for a solution to the
     $\cQ$-stabilization problem is given by Theorem \ref{T:1and2} (1).  As
     $\cQ$-stability always implies robust stability, the
     conditions in Theorem \ref{T:1and2} (1) give sufficient conditions for the
     existence of a static controller satisfying the robust
     stabilization problem.

     For the discussion of dynamic controllers some care must be taken,
     since the quantities $\delta_{1}, \dots, \delta_{d}$ are here
     uncertainties which are unknown to the controller.  To obtain
     sufficient conditions for the existence of a {\em  dynamic} controller
     solving the  robust stabilization problem, one only needs to apply
     the more flexible Theorem \ref{T:3} with the prescription that
     the controller state-space dimensions $n_{Kk}$ are to be equal to $0$
     for $k = 1, \dots, d$ but no constraint is imposed on $n_{KS}$
     (i.e., the controller is allowed to have dynamics corresponding to
     the frequency variable $\lambda$).  Similarly, the conditions in
     Theorem \ref{T:4} with the imposition that $n_{Kk} = 0$ for $k =
     1, \dots, d$ and $n_{KS}=0$ (the static controller case) or only
$n_{Kk}=0$ for
     $k = 1, \dots, d$ (the case where the controller is allowed to have
     dynamics with respect to the frequency variable) give sufficient
     conditions for the existence of a controller which solves the
     robust performance problem.

     As is now well-known (see \cite{Paganini, LZD, BGM3, DP}), if one
     expands the structured uncertainty  to include time-varying
     structured uncertainty, then robust stability is equivalent to
     $\cQ$-stability and the various conditions in Theorems
     \ref{T:1and2},
      and \ref{T:4} become necessary as well as sufficient for
     the existence of the respective type of controller solving the
     robust stabilization/performance problem.  The LFT model for this
     expanded uncertainty structure amounts to tensoring the system matrix
     \eqref{aggsys} with
     $I_{\ell^{2}}$ (the identity operator on the space $\ell^{2}$ of
     square-summable sequences)
     $$
     \begin{bmatrix} {\mathbf A} & {\mathbf B}_{1} & {\mathbf B}_{2} \\
         {\mathbf C}_{1} & {\mathbf D}_{11} & {\mathbf D}_{12} \\
         {\mathbf C}_{2} & {\mathbf D}_{21} & {\mathbf D}_{22}
     \end{bmatrix} =
     \begin{bmatrix} A & B_{1} & B_{2} \\ C_{1} & D_{11} & D_{12} \\
     C_{2} & D_{21} & D_{22} \end{bmatrix} \otimes I_{\ell^{2}}
     $$
     and expanding the block structure to have the form
      \begin{equation}  \label{tvDelta}
     \widetilde{\boldsymbol \Delta} = \left\{
     \widetilde \Delta = \begin{bmatrix} I_{n_{1}} \otimes \delta_{1}
     & & & \\
     & \ddots & & \\ & & I_{n_{d}} \otimes \delta_{d} & \\
     & & & I_{n} \otimes S \end{bmatrix} \colon
     \delta_{k} \in \cL(\ell^{2})  \text{ for } k = 1, \dots, d \right\}
     \end{equation}
     where $S$ is the shift operator on $\ell^{2}$
     $$  S \colon (c_{1}, c_{2}, c_{3}, \dots ) \mapsto (0, c_{1},
     c_{2} , \dots ).
     $$

     It can be shown that the results are unaffected if
     one replaces the shift operator $S$ in \eqref{tvDelta} by a general
     operator $\delta_{S} \in \cL(\ell^{2})$;  hence the LFT feedback
     model formalism carries over to this setting.

     \subsection{LPV systems}  \label{S:LPV}

     A second application of LFT models to robust stabilization and
     performance problems is in the context of gain-scheduling for
Linear-Parameter-Varying
     (LPV) systems.  We assume that we are given an LFT model of the
     form \eqref{aggsys} and \eqref{aggblock} where now the
     quantities $\delta_{1}, \dots, \delta_{d}$ are interpreted to be,
     rather than uncertainties, plant parameters varying
     in time. It is assumed that the controller  has access to these
     parameter values $\delta_{1}, \dots, \delta_{d}$ at each point in time
     $t$.  Then it makes sense to consider robust stabilization and
     robust performance problems where the controller is allowed to
     have dynamics in the uncertainty (now parameter) variables as well
     as in the frequency variable $\lambda$.  In this setting $\cQ$-stability is
     sufficient but not equivalent to robust stability.  We conclude
     that the conditions in Theorem \ref{T:1and2} (2) (adapted to the structure
     \eqref{aggsys} with \eqref{aggblock}) are sufficient for the
     existence of such a ``gain-scheduling'' controller (see
\cite{Packard}) which achieves
     robust stability, and, similarly, the conditions of Theorem
     \ref{T:4} (with constraints on the controller state-space
     dimensions $n_{K1}, \dots, n_{Kd}, n_{KS}$ at the discretion
     of the user) are sufficient for the existence of such a controller
     achieving robust performance.
     Theorem \ref{T:4} in this context is one of the main results of
     the paper \cite{AG}; the ``scaled-$H^{\infty}$ problem'' defined
     there is equivalent to finding a controller $\Sigma_{K}$ which
     achieves our ``$\cQ$-performance'' for the closed-loop
     system.

\end{document}